\documentclass[11pt]{amsart}
\usepackage{}
\usepackage{amssymb}
\theoremstyle{plain}
\newtheorem{Thm}{Theorem}

\newtheorem{Pro}[Thm]{Proposition}

\begin{document}

\title[Injectivity radius bound of Ricci flow ]
{ Injectivity radius bound of Ricci flow with positive Ricci
curvature and applications}

\author{Li Ma, Anqiang Zhu}
\address{Dr.Li Ma, Distinguished Professor, Department of mathematics \\
Henan Normal university \\
Xinxiang, 453007 \\
China} \email{nuslma@gmail.com}

\address{Anqiang Zhu, Department of mathematics \\Wuhan University,Wuhan 430072, China}
\email{anqiangzhu@yahoo.com.cn}

\thanks{The research is partially supported by the National Natural Science
Foundation of China 10631020 and SRFDP 20090002110019}

\begin{abstract}
In this paper, we study the injectivity radius bound for 3-d Ricci
flow. We also give an example showing that the positivity of Ricci
curvature can not be weakened into non-negativity. As applications
we show the long time existence of the Ricci flow with positive
Ricci curvature. We also partially settle a question in page 302
(line -13) of the book of Chow-Lu-Ni (2006).

{ \textbf{Mathematics Subject Classification 2000}: 53Cxx,35Jxx}

{ \textbf{Keywords}: injectivity radius bound, Ricci flow, positive
Ricci curvature, compactness}
\end{abstract}

 \maketitle

\section{Introduction}
  In this paper we study the injectivity radius bound for 3-d Ricci flow
  with positive Ricci curvature.
  We prove the following result.

\begin{Thm}\label{key}
Assume that $(M,g(t))$, $t\in [0,T)$, is a 3-d simply connected
complete non-compact n-dimensional Ricci flow with bounded
curvature,i.e. $|Rm(g)|\leq B$ for some constant $B>0$ on $M$ and
with positive Ricci curvature $Rc>0$. Then for each $t>0$,
$(M,g(t))$ is non-collapse everywhere on $M$. Namely, there is a
positive constant $c=c(t)>0$ such that the injectivity radius
$injrad(x,g(t))\geq c(t)$.
\end{Thm}
This result is inspired from the works \cite{M} \cite{MM} \cite{MC}.
  As applications, we continue our study of the global existence of the
  Ricci flow on the 3-dimensional complete
  non-compact Riemannian manifold $(M,g_0)$ with positive Ricci
  curvature. We assume that the Riemannian curvature of $(M,g_0)$  decays at
  infinity, namely,
  \begin{equation}\label{decay}
  |Rm(g_0)(x)|\to 0
  \end{equation}
  as the distance $d_{g_0}(x,x_0)\to \infty$.

We shall prove the following result.
\begin{Thm}\label{zhu}
Assume that $(M,g_0)$ is a 3-dimensional complete
  non-compact Riemannian manifold $(M,g_0)$ with positive Ricci
  curvature and with the condition (\ref{decay}). Then there is a
  global Ricci flow with the same properties of $(M,g_0)$.
\end{Thm}

It is known by results
  of R.Hamilton (\cite{Ham82} \cite{Ham95}) and W.Shi that the Ricci flow
  exists and preserve the above two conditions. R.Hamilton asks if
  the Ricci flow exists globally. The above result settles his
  question. Without the assumption of behavior about the curvature
  at infinity, there is a finite blow-up Ricci flow, which is the
  standard solution constructed by Perelman \cite{P}.

  We have proven the global existence of the Ricci flow under the positive curvature assumption on $(M,g_0)$.
The positive curvature condition is for the injectivity radius bound
used for the compactness theorem \cite{Ham95} \cite{L}.

We now recall the definition of type III Ricci flow. A global Ricci
flow $(M,g(t))$ is called Type III if there exists a constant $A>0$
such that
$$
\sup_{M^n\times [0,\infty)}t|Rm(g(t))|=A<\infty.
$$

We also partially settle a question in page 302 (line -13) of the
book \cite{Chow08}.

\begin{Thm}\label{ma}
Given a 3-d type III Ricci flow with positive Ricci curvature and
with a sequence $(x_i,t_i)\in M\times [0,\infty)$ and $ t_iK_i\geq c
$ where $K_i=|Rm(x_i,t_i)|$ for some positive constant $c>0$. Then
the family of Ricci flows
$$
g_i(t)=K_ig(t_i+K^{-1}_it)
$$
has a convergent subsequence in the sense of Cheeger-Gromov.
\end{Thm}

The interesting part of the result Theorem \ref{key} above lies in
that we don't assume any non-collapse for the initial metric. The
key part of our proof is to use the Gauss-Bonnet formula to the
minimal surface at infinity by using the simply connectedness at
infinity proved by Schoen-Yau \cite{SY}. We also point out that
there is a new interesting result of X.C. Rong \cite{rong} about
non-collapsing for Riemannian manifolds with almost-negative
curvature and with connectivity at infinity.

 Recall that for the metric $g$ with positive Ricci
curvature $Rc(g)>0$, the volume quotient
$$
\frac{Vol(B_g(x,,R)}{\omega_{n} R^n}
$$
is well defined monotone non-increasing function in $R\in
(0,\infty)$ via the Bishop theorem. One may use the lower bound of
this quotient to get the lower injectivity bound. However, near the
infinity, this quantity may be very small.

The key result Theorem \ref{key} is proved in next section. We also
give an example showing that the positivity of Ricci curvature can
not be weakened into non-negativity. The proofs of other results are
discussed in the last section.

\section{Proof of Theorem \ref{key} and its remarks}

Recall that for any $b>0$, the metric $g_b=b^2g$, we have the
following relations of the metric balls (see page 253 in
\cite{Chow08}) that
$$ B_{g}(x,r))=B_{g_b}(x,br), \ \
Vol_{g_b}(B_{g_b}(x,br))=b^nVol_g(B_{g}(x,r)),
$$
$$
|Rm{g_b}(x)|=b^{-2}|Rm{g}(x)|,
$$ and $$ injrad(x,g_b)=b\cdot injrad(x,g).
$$

We now prove Theorem \ref{key}.

\begin{proof}

We argue by contradiction. Assume the conclusion of Proposition
\ref{key} is not true for some $t_0\in (0,T)$. Then there exists
$x_j\to\infty$ such that
$$
injrad(x_j,g(t_0))\to 0
$$
and a minimizing geodesic $\sigma_j$ based at $x_j$ such that
$L(\sigma_j)=2injrad(x_j,g(t_0))$. Set
$\lambda_j=injrad(x_j,g(t_0))$. We now normalize the Ricci flow
$g(t)$ at $(x_j,t_0)$ by
$$
g_j(t)=\lambda_j^{-2}g(\lambda_j^2t+t_0).
$$
Then the Ricci flow $g_j(t)$ has
$$injrad(x_j,g_j(0))=1$$ and
$$
|Rm(g_j(t))|\leq B\lambda_j^2\to 0.
$$
Using the Cheeger-Gromov-Hamilton convergence theorem we may assume
that
$$
(M,g_j(t),x_j)\to (M_\infty, g_\infty(t),x_\infty)
$$
in the sense of $C^2$ Cheeger-Gromov sense
(\cite{Ham95},\cite{Cho06}, \cite{Chow08}) and $Rm(g_\infty)=0$.
Hence $B_{g_\infty(0)}(x_\infty,1)$ is the unit euclidean ball and
$$Vol(B_{g_\infty(0)}(x_\infty,1))=\omega_{n}.$$

We now recall a key fact from Schoen-Yau's work \cite{SY} about the
geometry of the complete non-compact Riemannian manifold $(M^3,g)$
of positive Ricci curvature that $M$ is diffeomorphic to $R^3$.
Recall here that $M$ is simply connected at infinity if there are no
compact set $K\subset M$ and a sequence of Jordan curves
$\{\sigma_j\}$ tending uniformly to infinity such that any sequence
of disks $\{D_j\}$ with $\partial D_j=\sigma_j$ has $D_j\bigcap
K\neq\phi$ for each $j$.

Using $injrad(x_j,g_j)=1$ we get a closed minimizing geodesic
$\gamma_j=(\sigma_j)$ in $(M,g_j)$ based at $x_j$ (and the limit
$\gamma_\infty$ of $\gamma_j$ is a true closed geodesic in flat
space $(M_\infty,g_\infty)$). Smoothing $\gamma_j$ (at $x_j$, which
may be a corner point) and using the simply-connectedness at
infinity of $(M,g(t_0))$ (by the well-known theorem of Schoen-Yau
\cite{SY}) we can bound it by a minimizing (immersed) minimal disk
$\Sigma_j$ (and its existence is guaranteed by the works of
S.Hildebrandt and C.B.Morrey) with its area bounded above by some
uniform constant, saying $8\pi$ for large $j$. We remark that the
area bound comes from a reversible construction from the fact (Lemma
9.5.1(b) in \cite{Morrey}) that $\lim_{j\to\infty}
A(\Sigma_j)=A(\Sigma_\infty)$, where $\Sigma_\infty$ is the
minimizing minimal disk bounded by the loop $\gamma_\infty$ in
$(M_\infty,g_\infty)$. One can also get the area bound by using the
pull-back of the surface $\Sigma_\infty$ to the space $(M,g_j(t))$
as the comparison surface. The geometric picture of this
construction is that the limit space $M_\infty$ can not be
$(S^1\times R)\times R$. For otherwise, we would have that
$\gamma_j$ has its limit $S^1$ and the $\Sigma_j$ converges to
$S^1\times R$. Hence, the area of $\Sigma_j$ would be large. Since
the limit of Busemann functions $B_{\gamma_j}$ corresponding to
$\gamma_j$ is the coordinate function in the real line $R$, we can
use the cut of regular level sets of the Busemann function
$B_{\gamma_j}$ to strictly reduce $\Sigma_j$ in its area, which is a
contradiction to the minimizing property of $\Sigma_j$ (being an
area minimizing disk). See also the remark below and the references
\cite{Ham99} \cite{P} for related idea. On the surface $\Sigma_j$
with the unit normal vector $N$ and the second fundamental form
$A_j$ of $\Sigma_j\subset M$, we know that
$$
Rc_{g_j}(N,N)+\frac{1}{2}|A_j|^2=\frac{R(g_j)}{2}-K_j.
$$
Here $K_j$ is the Gauss curvature of the minimal surface $\Sigma_j$.
Then using $Rc>0$ we have
$$
\frac{1}{2}|A|^2+K_j\leq\frac{R(g_j)}{2}\to 0
$$
and by using the uniform area bound of $\Sigma_j$,
$$
\int_{\Sigma_j}K_j\leq \frac{1}{2}\int_{\Sigma_j}R(g_j)\leq
\sup_{\Sigma_j} R(g_j) A(\Sigma_j)\to 0.
$$

 Recall that the Gauss-bonnet formula (and here we may assume that $\Sigma_j$'s are embedded)
 $$
\int_{\Sigma_j}K_j+\int_{\gamma_j}k_{g_j}=2\pi,
 $$
 which gives us a contradiction since both terms in left side tend
 to zero or less than  zero as $j\to\infty$.

We then complete the proof of Theorem \ref{key}.
\end{proof}

We give two remarks about the argument above. One is that the
important step in getting the minimizing minimal disk bounding
$\gamma_i$ with uniform bounded area in the above argument is using
the simply-connectedness at infinity of $M$ and non-splitting of
$M$. Recall here that for a given ray $\widetilde{\gamma}(t)$ with
$\widetilde{\gamma}(0)=x_0$, the Busemann function is defined by
$$
B_{\widetilde{\gamma}}(x)=\lim_{t\to\infty}
[t-d(x,\widetilde{\gamma}(t))],
$$
where $d(x,x_0)$ is the distance function in $(M,g)$.  The sequence
$(x_i)$ determines a ray $\widetilde{\gamma}(t)$ and we write by
$B_{\widetilde{\gamma}}$ ($B_i$) the corresponding Busemann function
in $(M,g)$ (in $(M,g_i,x_i)$). By the simply-connectedness at
infinity (based on the fact that the level sets of the Busemann
function $B_i$ are of positive mean curvature, see line 7 in page
221 in \cite{SY}, and correspondingly the level sets of the distance
function $d(x):=d(x,x_0)$ are of almost negative mean curvature when
$d(x)$ large) and Bernstein type theorem (Theorem 2 in \cite{SY}) we
can always find minimizing minimal disks $D_i$ bounding $\gamma_i$
in $(M,x_i)$ (according to the works of Hildebrandt \cite{Hi} and
C.B.Morrey \cite{Mo}) such that there are only finite $D_i$ can
interest any fixed compact subset $K\subset M$ (and otherwise there
is a non-trivial stable minimal surface in $(M,g)$, which is
impossible \cite{SY}). We can show that for
$dist(\gamma_i,\gamma_j)$ large, there is a contractible domain
$\Omega(i,j)$ such that both $D_i$ and $D_j$ are in the boundary of
$\Omega(i,j)$. This implies that we find at least one comparison
disk $G_i$ spanning $\gamma_i$ with uniform area bound such that it
lies between two level sets $d^{-1}(a_i)$ and $B_i^{-1}(a_i+L_i)$
for some $a_i>>1$ and $0<L_i\leq 2$ (because the half length of
$\gamma_i$ is one). Another way to find such comparison surface is
below. Note that the limit surface of $D_i$ is a flat disk in the
flat space $(M_\infty,g_\infty,x_\infty)$. Using this limit surface
we can also construct a comparison surface to each $D_i$ such that
the area of $D_i$ is uniformly bounded. This step can not be carried
through to 3-d manifold $M$ with non-negative Ricci curvature since
we may not find a minimizing disk bounding $\gamma_i$ in $M_i$ with
uniform bounded area. Here is the example that $M=N\times R$ with
$N=R^2$ equipped with cigar metric $h$. With this metric, $M$ has
non-negative Ricci curvature and it does not have the strong
contractible property. In this example, we let $x_i$ in N going to
infinity. Let $\gamma_i$ be a closed geodesic realizing $injrad(x_i)
\approx 1$. The area of the disk $D_i$ bounding $\gamma_i$ goes to
infinity with $i\to\infty$. Moreover $(N,h,x_i)$ converges to a flat
cylinder $S^1\times R$ and $\gamma_i$ converges to the $S^1$ factor
of the cylinder, which does not bound a disk. Now if $M = R^3 =
N\times R$ and  $g=h+ds^2$ on M, consider $(0,x_i)$ in M going to
infinity and the limit here is the flat $ S^1\times R^2$, and
$\gamma_i$ converges to the $S^1$ factor, which does not bound a
disk. The key point for this example is that with the metric
$g=h+ds^2$, $R^3$ has a splitting structure, which makes the
minimizing disk have no area bound.

The other is that once we have the uniform injectivity radius bound
of $g(t_0)$, for some $t_0>0$, we have the the uniform injectivity
radius bound $c(t_0, B)$ of $g(t)$ for $t>t_0$. The constant
$c(t_0,B)$ depends only on $t_0,B$ by using the estimates of Carron
Proposition 4 in \cite{C}).

\section{Remarks about proofs of other results}

The key step in the argument about the generalized version of the
result in \cite{MZ} is about the injectivity radius bound, which now
can be obtained via Theorem \ref{key}. Then the proof of Theorem
\ref{zhu} is almost the same as in \cite{MZ}. So we omit the detail.

The proof of Theorem \ref{ma} follows easily from Hamilton's
compactness theorem of Ricci flow \cite{Ham95} \cite{Chow08}.


\begin{thebibliography}{20}
\bibitem{BBB}
Bessieres, L.; Besson, G.; Boileau, M.; Maillot, S.; Porti, J.
\emph{Collapsing irreducible 3-manifolds with nontrivial fundamental
group}. Invent. Math. 179 (2010), no. 2, 435-460


\bibitem{C}
G. Carron, \emph{Inegalities isoperimetriques de Faber-Krahn et
consequences}, Semin. Congr., 1, Soc. Math. (France, Paris, 1996)
205-232

\bibitem{Che82}
J.Cheeger, M.Gromov, M.Taylor, \emph{Finite propagation speed,
kernel estimates for functions of the Laplace operator, and the
geometry of complete Riemannian manifolds}. J.Diffential.Geom.
17(1982), 15-53


\bibitem{Cho06}
 B.Chow, P.Lu, L.Ni, \emph{Hamilton's Ricci Flow}. Science
Press. American Mathematical Society, Beijing,Providence, 2006.

\bibitem{Chow08}
B.Chow, etc,\emph{The Ricci flow: The techniques and Applications,
Part I: geometric aspects}, American Mathematical Society,
Beijing,Providence, 2008.



\bibitem{Ham82}
R.Hamilton, \emph{Three-manifolds with positive Ricci curvature}.
 J. Differential Geom., 2(1982)255-306.

 \bibitem{Ham99}
 R.Hamilton,
\emph{Non-singular solutions of the Ricci flow on three-manifolds}.
Comm. Anal. Geom. 7 (1999), no. 4, 695-729.


\bibitem{Ham95}
R. Hamilton, \emph{The formation of singularities in the Ricci
flow}, Surveys in Differential Geometry, Vol II, International Press
(1995), 7-136.

\bibitem{Hi} S.Hildebrandt, \emph{Boundary behavior of minimal surfaces}. Arch. Rat.
Mech. Anal. 35 (1969), 47-82.


\bibitem{H}
H.Huang, \emph{A note on Ricci flow on non-compact manifolds},
Journal of Math. Study, 42(4)(2009)351-356.

\bibitem{L}
J.Lott, \emph{On the long time behavior of type III Ricci flow
solutions}, Math. Ann. 339(2007), 627-666.


\bibitem{M}
L.Ma, A complete proof of Hamilton's conjecture,
http://arxiv.org/abs/1008.1576v1


\bibitem{MM}
L.Ma, \emph{Expanding Ricci solitons with pinched Ricci curvature},
Kodai Math.Journal, 34 (2011), 140-143

\bibitem{MC} L.Ma, L.Cheng,\emph{
Yamabe flow and Myers type theorem on complete manifolds},
J.Geom.Anal.,DOI 10.1007/s12220-012-9336-y, online, 2012

\bibitem{MZ}
L.Ma, A.Zhu,  \emph{Nonsingular Ricci flow on a noncompact manifold
in dimension three}, C.R.Mathematique, ser.I,137(2009)185-190.

\bibitem{Mo} C.B.Morrey, \emph{The problem of Plateau on a riemannian manifold}. Ann.
Math. 49 (1948), 807-851.

\bibitem{Morrey}
C.B.Morrey,Jr., \emph{Multiple Integrals in the Calculus of
Variations}, Grundlehren Math. Wiss. 130, Springer-verlag, Berlin,
1966 edition.

\bibitem{P}
G. Perelman, \emph{Finite extinction time for the solutions to the
Ricci flow ow on certain three-manifolds}. arXiv:math.DG/0307245.

\bibitem{rong}
X.C.Rong, \emph{Almost non-negative curvature vs. collapse in
dimension three}, personal commnication,2011.

\bibitem{SY}
R. Schoen, S.T.Yau, \emph{Complete three-dimensional manifolds with
pos- itive Ricci curvature and scalar curvature}. Seminar on
Differential Geometry, pp. 209-228, Ann. of Math. Stud., 102,
Princeton Univ. Press, Princeton, N.J., 1982.

\bibitem{Shi89b}
W. X. Shi,\emph{Ricci deformation of metric on complete noncompact
Riemannian manifolds}, 30(1989)303-394; \emph{Complete noncompact
three manifolds with nonnegative Ricci curvature}, JDG,
29(1989)353-360.


\end{thebibliography}
\end{document}